\documentclass[11pt]{article}
\usepackage{latexsym}
\usepackage{amsfonts}
\usepackage{epsfig}
\usepackage{amssymb}
\usepackage{amscd}
\usepackage[centertags]{amsmath}
\usepackage[all]{xy}
\usepackage{amsthm}
\usepackage{psfrag}

 \newtheorem{teo}{Theorem}[section]
 \newtheorem{cor}[teo]{Corollary}
 \newtheorem{lem}[teo]{Lemma}
 \newtheorem{pro}[teo]{Proposition}
 \theoremstyle{definition}
 
 \theoremstyle{remark}
 \newtheorem{rem}[teo]{Remark}
 \newtheorem{exam}[teo]{Example}
 \numberwithin{equation}{section}

\newcommand{\m}{\mathfrak{m}}
\newcommand{\D}{\mathcal{D}}
\newcommand{\K}{\mathcal{K}}
\newcommand{\C}{\mathcal{C}}
\newcommand{\E}{\mathcal{E}}
\newcommand{\G}{\Gamma}
\newcommand{\OO}{\mathcal{O}}
\newcommand{\p}{\textsl{p}}
\newcommand{\q}{\textsl{q}}
\renewcommand{\u}{\textsl{u}}
\renewcommand{\v}{\textsl{v}}
\newcommand{\w}{\textsl{w}}
\newcommand{\I}{\mathbf{I}_R}

\newcommand{\nop}{}

\newif\ifprivate
\privatetrue

\def\???{\ifprivate {\bf {???}} \marginpar{{\Huge {\bf ?}}}
\else \fi}

\begin{document}

\title{Equisingularity classes of birational projections of normal singularities to a plane}

\author{Maria Alberich-Carrami\~nana and Jes\'{u}s Fern\'{a}ndez-S\'{a}nchez
\thanks{This research has been partially supported by the Spanish
Committee for Science and Technology (DGYCIT), projects
MTM2005-01518 and MTM2006-14234-C02-02, and the Catalan Research
Commission. The second author completed this work as researcher of
the program \emph{Juan de la Cierva}.
} \\
\small{Departament de Matem\`atica Aplicada I, Universitat Polit\`ecnica de Catalunya,}\\
\small{Av. Diagonal 647, 08028-Barcelona, Spain.} \\ \small{e-mail: maria.alberich@upc.edu, jesus.fernandez.sanchez@upc.edu}
    }


\maketitle
%
\begin{abstract}
Given a birational normal extension $\OO$ of a two-dimensional local regular ring $(R,
\m)$, we describe all the equisingularity types of the complete $\m$-primary ideals $J$
in $R$ whose blowing-up $X=Bl_J(R)$ has some point $Q$ whose local ring $\OO_{X,Q}$ is
analytically isomorphic to $\OO$.
\end{abstract}

\section*{Introduction}
A sandwiched surface singularity $(X, Q)$ is a normal surface
singularity that can be projected birationally to a non-singular
surface. From a more algebraic point of view, the local ring $\OO$
of any sandwiched singularity is a birational normal extension of
a two-dimensional local regular ring $R$.
Once a sandwiched surface singularity has been fixed, in this
paper we address the problem of describing the equisingularity
classes of all its birational projections to a plane.
The problem of classifying the germs of sandwiched surface
singularities was already posed by Spivakovsky. As he claims in
\cite{Spivak} this problem has two parts: discrete and continuous.
The continuous part is to some extent equivalent to the problem of
the moduli of plane curve singularities, while the main result of
this paper solves completely the combinatorial part.
Any birational projection from a sandwiched singularity to a
plane is obtained by the morphism of blowing up a complete
$\m_O$-primary ideal in the local ring of a regular point $O$ on the plane.
Our goal is to give all the equisingularity types of these ideals.
Namely, fixed a birational normal extension $\OO$ of a local
regular ring $(R, \m_O)$, we describe the equisingularity type of
any complete $\m_O$-primary ideal $J\subset R$ such that its
blowing-up $X=Bl_J(R)$ has some point $Q$ whose local ring
$\OO_{X,Q}$ is analytically isomorphic to $\OO$. In this case, we
will say that the surface $X$ contains the singularity $\OO$ for
short, making a slight abuse of language.
This is done by describing the Enriques diagram of the cluster of
base points of any such ideal $J$: such a diagram will be called
an Enriques diagram for the singularity $\OO$.
Recall that an Enriques diagram is a tree together with a binary
relation (proximity) representing the topological equivalence
classes of clusters of points in the plane (see \S
\ref{ED-dualG}).
Previous works by Spivakovsky \cite{Spivak} and M\"{o}hring
\cite{Mohring} describe a type of Enriques diagram that exists for
any given sandwiched surface singularity (detailed in \S 2) and
provide other types mostly in the case of cyclic quotients (see
\cite{Mohring} 2.7) and minimal singularities (see \cite{Mohring}
2.5).

The organization of the paper is as follows. Section 1 is devoted
to recalling some definitions concerning the language of infinitely
near points, sandwiched surface singularities and graphs.
Fundamental for our purpose will be the notion of Enriques
diagram, introduced in \cite{EC85}. In Section 2, after some
technical results, we introduce the concept of \emph{contraction}
for a sandwiched surface singularity $\OO$. By a contraction we
mean the resolution graph $\G_{\OO}$ of $\OO$ (a sandwiched graph,
as introduced in \cite{Spivak}) enriched by some proximities
between their vertices, these proximities being compatible with
the weights of the graph. Fixed a sandwiched graph, the problem of
finding the whole list of possibilities for such proximities is
the hard part of our work. This is achieved in Section 3, by
proving that any contraction for $\G_{\OO}$ may be recovered from
some contraction of the graph obtained from $\G_{\OO}$ by removing
one end. This fact is the key result in order to describe a
procedure to obtain all the contractions for $\OO$. Finally, in
Section 4, we explain how to complete contractions in order to
obtain any Enriques diagram for $\OO$.

\vspace{3mm}
\noindent \textbf{Acknowledgement}
The authors thank E. Casas-Alvero for drawing their attention to the problem addressed in this paper. They also thank M. Spivakovsky for the conversations held on this topic.

\section{Preliminaries}

In this section, we fix notation and recall some of the facts
concerning sandwiched surface singularities and base points of
ideals which will be used throughout this paper, and we focus on
our problem. A standard reference for most of the tools and
techniques treated here is the book by Casas-Alvero \cite{Cas00}.

\subsection{Infinitely near points and complete ideals}

Let $(R,{\m_O})$ be a regular local two-dimensional
$\mathbb{C}$-algebra and $S=Spec(R)$.
A \emph{cluster} of points of $S$ with origin $O$ is a finite set
$K$ of points infinitely near or equal to $O$ such that, for any
$p\in K$, $K$ contains all points to which $p$ is infinitely near.
A subset of a cluster $K$ is a \emph{subcluster} if it is a
cluster (with origin some point of $K$).
By assigning integral multiplicities $\nu=\{\nu_p\}$ to the points
of $K$, we obtain a \emph{weighted cluster} ${\mathcal
K}=(K,\nu)$; the multiplicities $\nu$ are called the
\emph{virtual multiplicities} of $\K$.
A point $p$ is said to be \emph{proximate} to another point $q$ if
$p$ is infinitely near to $q$ and lies on the strict transform of
the exceptional divisor of blowing up $q$.
We write $p\geq q$ if
$p$ is infinitely near or equal to $q$, and $p\rightarrow q$ if
$p$ is proximate to $q$.
The relation $\geq$ is an ordering of the infinitely near points,
and it will be considered as their natural ordering.
A point $p\in K$ is \emph{free} if it is proximate to only one
point, which is necessarily the immediate predecessor, and $p$ is
\emph{satellite} if it is proximate to two points; otherwise, the
point is necessarily the origin of the cluster.
The number $\rho^{\K}_p=\nu_p- \sum_{q \rightarrow p} {\nu _q }$
is the \emph{excess} at $p$ of ${\mathcal K}$. \emph{Consistent
clusters} are those weighted clusters with non-negative excesses
at all their points. We write $\K_+$ for the set of the
\emph{dicritical} points of $\K$, that is, the points with
positive excess.

If $\mathcal{K}=(K, \nu)$ and $\mathcal{K}'=(K', \nu ')$, define
the \emph{sum} $\mathcal{K}+ \mathcal{K}'$ as the weighted cluster
whose set of points is $K \cup K'$ and whose virtual
multiplicities are $\nu_{p}+\nu'_{p}$ for $p \in K \cup K'$
(\cite{Cas00} 8.4). This operation is clearly associative and
commutative, thus making the set of all weighted clusters with
origin at $O$ a semigroup.
Consider the set $\mathbf{W}$ of all consistent clusters with
origin at $O$ with positive virtual multiplicities.
Again, $\mathbf{W}$ equipped with the sum is clearly a semigroup.
A weighted cluster $\mathcal{K} \in \mathbf{W}$ is called
\emph{irreducible} if it is so as element of the semigroup
$\mathbf{W}$, that is, $\mathcal{K}$ is not the sum of two
elements of $\mathbf{W}$.
To any point $p$, $p \ge O$, we associate the irreducible cluster
$\mathcal{K}(p)$ in $\mathbf{W}$ which has virtual multiplicity
one at $p$, which will be called the \emph{irreducible cluster in}
$\mathbf{W}$ \emph{ending at} $p$.

Two \emph{clusters} $K$ and $K'$ are called \emph{similar} if
there is a bijection (\emph{similarity}) $\varphi: K
\longrightarrow K'$ so that both $\varphi$ and $\varphi ^{-1}$
preserve ordering and proximity.
Two \emph{weighted clusters} $\mathcal{K}=(K, \nu)$ and
$\mathcal{K}'=(K', \nu ')$  are called \emph{similar} if there is
a similarity between $K$ and $K'$ preserving virtual
multiplicities (see \cite{Cas00} 8.3).
An analytic isomorphism $\Phi$ defined in a neighborhood of $O$ clearly induces a
similarity between each cluster $K$ with origin $O$ and its image $\Phi(K)$ (\cite{Cas00}
3.3).
Furthermore, if $\Phi$ is only a homeomorphism, then $K$ and $\Phi(K)$ are still similar
(\cite{Cas00} 8.3.12)

If $\pi_K:S_K\longrightarrow S$ is the composition of the blowing-ups of all points in
$K$, write $E_K$ for the exceptional divisor of $\pi_K$ and $\{E_p\}_{p\in K}$ for its
irreducible components. We denote by $\mathbf{A}_K=(E_p\cdot E_q)_{p,q\in K}$ the
\emph{intersection matrix} of $E_K$: if $p=q$, its coefficient is just the
self-intersection of $E_p$, and equals $-r_p-1$, where $r_p$ is the number of points in
$K$ proximate to $p$; if $p\neq q$, $E_p\cdot E_q=1$ in case $E_p\cap E_q\neq \emptyset$,
and $E_p\cdot E_q=0$ otherwise. It can be easily seen that $E_p\cap E_q\neq \emptyset$ if
and only if $p$ is maximal among the points of $K$ proximate to $q$ or vice-versa (cf.
\cite{Cas00} 4.4.2).
Notice that $\mathbf{A}_K$ is an invariant of the similarity class of $K$.

If ${\mathcal K}$ is a weighted cluster, there is a well established notion for a germ of
curve to go through $\mathcal{K}$ (which is a linear condition, see \cite{Cas00} 4.1),
and the equations of all curves going through $\mathcal{K}$ define a complete
$\m_O$-primary ideal $H_{\mathcal K}$ in $R$ (see \cite{Cas00} 8.3). Any complete
$\m_O$-primary ideal $J$ in $R$ has a \emph{weighted cluster of base points}, denoted by
$BP(J)$, which consists of the points shared by, and the multiplicities of, the curves
defined by generic elements of $J$. Moreover, the maps $J\mapsto BP(J)$ and ${\mathcal
K}\mapsto H_{\mathcal K}$ are reciprocal isomorphisms between the semigroup
$\textbf{I}_R$ of complete $\m_O$-primary ideals in $R$ (equipped with the product of
ideals) and the semigroup $\textbf{W}$ (see \cite{Cas00} 8.4.11 for details).
If $p \ge O$, denote by $J(p)$ the ideal in $\mathbf{I}_R$
corresponding by the preceding isomorphism to the irreducible
cluster $\mathcal{K}(p) \in \mathbf{W}$ ending at $p$, that is,
$J(p) = H_{\mathcal{K}(p)}$.

A couple of ideals $J, J'$ in $\mathbf{I}_R$ are \emph{equisingular} if $BP(J)$ and
$BP(J')$ are similar (\cite{Cas00} 8.3). Notice that two equisingular complete ideals in
$\mathbf{I}_R$ have equisingular (that is, topologically equivalent) generic germs and
equal codimensions (\cite{Cas00} 8.3.9).


\subsection{Sandwiched surface singularities} \label{section 1.1}

The main references here are \cite{Spivak} and \cite{FS1}.
If $I\in \I$, we denote by $\pi_I:~X=Bl_I(R)\longrightarrow S$ the blowing-up of
$I$. The surface $X$ is not regular in general, and its singularities are sandwiched
singularities.
Moreover, if $K$ is the set of base points of $I$, we have a
commutative diagram
\begin{equation}\label{CommDiagr}
  \xymatrix{  {S_K} \ar[r]^{f}\ar[rd]_{\pi_K} & {X} \ar[d]^{\pi_I} \\ &  {S}}
\end{equation}
where the morphism $f$, given by the universal property of the blowing-up, is the minimal
resolution of the singularities of $X$ (\cite{Spivak} Remark 1.4).
Let $\OO$ be any singularity of $X$; then we say that $I$ is an
\emph{ideal for }$\OO$.
It follows that the exceptional divisor $E_{\OO}$ associated with
the minimal resolution of $\OO$ is a connected subset of the
exceptional divisor $E_K$.
There is a bijection between the set of irreducible components of
$\pi_I^{-1}(O)$ and the set of dicritical points of $\K =BP(I)$
(see \cite{FS4,Lip69}). This allows to write $\{L_p\}_{p\in \K_+}$
for the set of these components on $X$. Because of this, we may
think of $\OO$ as a singularity obtained by contracting a
connected curve (which will be called $E_{\OO}$) of $E_K$
containing no component with self-intersection $-1$ (such a
component $E_p$ is necessarily the exceptional divisor of the
blowing-up of some maximal point of $K$ and thus, a dicritical
point).

For any ideal $J= \prod_{p \in \K_+}  J(p)^{\alpha(p)}$ with
positive $\alpha(p)$, we have an analytic isomorphism $X \cong
Bl_J(R)$ (cf. \cite{Spivak}, Corollary I.1.5).
Since we are interested in sandwiched singularities modulo
analytic isomorphism, the relevant information we need to retain
about $\K=BP(I)$ is, on one hand, its set of points $K$ and, on
the other, knowing which of the points of $K$ are dicritical
(the rest being non-dicritical, of excess zero).

\subsection{Enriques diagrams and dual graphs} \label{ED-dualG}

We introduce the Enriques diagrams and the weighted dual graphs
related to them.
The Enriques diagrams are combinatorial objects that enclose the
topological information of the clusters of infinitely near points
in $S$, namely they represent the similarity classes of clusters.

A \emph{tree} is a finite graph with a partial order relation $\leq$ between the
vertices, without loops, which has a single initial vertex, or \emph{root}, and every
other vertex has a unique immediate predecessor. The vertex $q$ is said to be a
\emph{successor} of $p$ if $p$ is the immediate predecessor of $q$. If $p$ has no
successors then it is an \emph{extremal vertex}.
The set of vertices of a graph will be denoted by the same letter
as the graph itself.
An \emph{Enriques diagram} $D$ (\cite{EC85} Enriques IV.I, \cite{Cas00} Casas 3.9; see
also \cite{GSG92} and \cite{KP99} for a combinatorial presentation) is a tree with a
binary relation between vertices, called \emph{proximity} and denoted by
$\rightarrow_{D}$, which satisfies:

\begin{enumerate}
 \item Every vertex but the root is proximate to its immediate predecessor; the root is
 proximate to no vertex.

 \item If $p\rightarrow_D q$, then $p>q$ and there is at most one other vertex in $D$
 proximate to both of them.

 \item Any vertex is proximate to at most two other vertices.
The vertices which are proximate to two points are called
\emph{satellite}, the other vertices, but the root, are called
\emph{free}. If $q$ is the immediate predecessor of $p$, and
$p\rightarrow_D q'$, then $q\rightarrow_D q'$.

\end{enumerate}

If $p$ is a vertex in $D$, we write $r_D(p)$ for the number of vertices in $D$ proximate to $p$.
A satellite vertex is said to be \emph{satellite of} the last free
vertex that precedes it.
In order to express graphically the proximity relation, Enriques
diagrams are drawn according to the following rules:
\begin{enumerate}
\item If $q$ is a free successor of $p$
   then the edge going from $p$ to $q$ is smooth and curved and,
   if $p$ is not the root, it has at $p$ the same tangent
   as the edge joining $p$ to its predecessor.
\item The sequence of edges connecting a maximal
   succession of vertices proximate to the same vertex $p$
   are shaped into a line segment, orthogonal to the edge joining $p$
   to the first vertex of the sequence.
\end{enumerate}

If $K$ is a cluster, there is an Enriques diagram $D_K$ naturally
associated with it by taking one vertex for each point of $K$ and
the proximity of the cluster as the proximity of $D_K$; conversely,
for any Enriques diagram $D$ there is some cluster $K$ with origin
$O$ whose Enriques diagram $D_K$ is $D$. If no confusion may
arise, we will label the points in $K$ and their corresponding
vertices in $D_K$ with the same symbol.
A connected subtree of an Enriques diagram $D$ is a
\emph{subdiagram} if it is an Enriques diagram with root some
vertex of $D$ and whose proximity is the restriction of the
proximity of $D$.
Observe that $K'$ is a subcluster of $K$ if and only if the associated Enriques
diagram $D_{K'}$ is a subdiagram of $D_K$.
If $D$ is the Enriques diagram associated with $K$ and $p \in K$,
we denote by $D(p)$ the Enriques diagram of the irreducible
cluster $\mathcal{K}(p)$ ending at $p$.
If $p$ is extremal, $D(p)$ is called a \emph{branch} of $D$.

By assigning to an Enriques diagram $D$ a \emph{marking map}
$\rho: D\rightarrow \{ + , 0 \}$, we obtain a \emph{marked
Enriques diagram} $\D=(D,\rho)$. Any consistent cluster $\K$
induces a marking map $\rho: D_K\rightarrow \{ + , 0 \}$ by taking
$\rho(p)=+$ if $p$ corresponds to a dicritical point of $\K$ (in
this case, $p$ is called a \emph{dicritical vertex}), and
$\rho(p)=0$ otherwise. A \emph{marked subdiagram} $\D' = (D',
\rho') $ of $\D$ is a marked Enriques diagram where $D'$ is a
subdiagram of $D$ and $\rho'$ is the restriction of $\rho$ to
$D'$. Observe that the extremal vertices of a marked Enriques
diagram associated with some $\K\in \mathbf{W}$ are always
dicritical.

If $\OO$ is a sandwiched surface singularity, we say that
$\mathcal{D}$ is an \emph{Enriques diagram for} $\OO$ if it is the
marked Enriques diagram of $BP(I)$, for some ideal $I\in \I$ for
$\OO$.
Under this framework, the goal of this paper is to describe all
the Enriques diagrams for a given $\OO$.

Incidence between the irreducible components of a divisor $E$ on a
surface is usually represented by means of the \emph{weighted dual
graph} of $E$. It is defined by taking a vertex for each component
of $E$, and by joining two vertices by an edge if and only if the
corresponding components of $E$ meet; each vertex is weighted by
taking minus the self-intersection of the corresponding component.
If $D$ is the Enriques diagram of a cluster $K$, the \emph{(weighted)
dual graph} of $D$, denoted by $\Gamma _D$, is the weighted dual
graph of the exceptional divisor $E_K$ (which has no loops).
Since the information enclosed in the weighted dual graph
is the same as that contained in the intersection matrix of $K$, this definition is consistent.
\begin{rem}
The similarity class of a cluster may be represented
either by its Enriques diagram or by its weighted dual graph,
since from the intersection matrix the ordering $\le$ (of being
infinitely near) and the proximity may be inferred. In fact, this is also true for rational surface singularities. From the intersection matrix $\mathbf{A}$ of a rational surface singularity, the fundamental cycle $Z$ may be computed (see \cite{Lau72} Theorem 4.2) and from it, the order of the blowing-ups performed to resolve the singularity: the negative entries of $\mathbf{A}Z$ correspond to the exceptional components having appeared in the last blowing-up (cf. Theorem 1.14 of \cite{Reguera97}).
It is worth noticing that the proximity of $D$ cannot be recovered
in general only from its dual graph without weights (see
\cite{Cas00} 4.4).
\end{rem}
A \emph{non-singular graph} is the weighted dual graph of some
Enriques diagram (cf. \cite{Spivak}).
The vertex in $\Gamma_{D}$ corresponding to $p$ in $D$
will be denoted by $\p$, written in \textsl{roman} font.
By the (weighted) dual graph of a marked Enriques diagram
$\D=(D,\rho )$ we mean the dual graph of $D$ and it will
also be denoted by $\G_D$. The vertices of $\G_D$ corresponding to
dicritical vertices (non-dicritical, respectively) of $\D$ will be
called \emph{dicritical} (\emph{non-dicritical}, respectively),
too.
If $I$ is a complete $\m_O$-primary ideal in $R$, we will write
$\D_I$ and $\Gamma_I$ to mean the marked Enriques diagram and the
weighted dual graph of its cluster of base points $BP(I)$,
respectively.

\begin{rem}\label{w=1extrem}
The weighted dual graph $\G_{D}$ can be constructed as follows:
take one vertex in $\Gamma_{D}$ for each vertex of $D$, and
connect two vertices in $\Gamma_{D}$ by an edge if and only if one
of the corresponding vertices in $D$ is maximal among the vertices
in $D$ proximate to the other. Moreover, if $\p$ is a vertex of
$\Gamma_{D}$, its weight $\omega(\p)$ is $r_D(p) + 1$ (cf. \S 4.4
of \cite{Cas00} for details).
A vertex $\p \in \Gamma_{D}$ has weight $\omega(\p)=1$ if and only if $\p$ is
extremal in $D$.
\end{rem}
A \emph{chain} $ch_{\G}(\q,\p)$ of a graph $\G$ without loops is
the subgraph composed of all vertices and edges between the
vertices $\q,\p\in \G$; it will be described by the ordered
sequence of vertices between $\q$ and $\p$, and $d_{\G}(\q,\p)$
will denote its length. Two vertices $\q,\p\in \G$ are
\emph{adjacent} if $d_{\G}(\p,\q)=1$; a vertex is an \emph{end} if
it is adjacent to only one vertex.
A \emph{weighted subgraph} of a weighted graph $\G$ is a subgraph of $\G$ whose
vertices have the same weights as $\G$.

The following result describes the proximity relations between the vertices of a chain:

\begin{lem}
\label{chain} Let $q \leq p$ be two vertices of an Enriques
diagram $D$, and consider the non-singular graph $\G$ of $D$.
\begin{itemize}
\item[(a)] If $\u\in ch_{\G}(\q,\p)$, then $q \leq u$; if $u \neq
p$, either $u \leq p$ or $p \leq u$. Moreover, all the vertices of
$ch_{\G}(\q,\p)$ correspond to vertices in the same branch of $D$.
\item[(b)] Write $ch_{\G}(\q,\p)=\{\u_0=\q,\u_1,\ldots,\u_n,\u_{n+1}=\p\}$. There exists
some $i_0\in \{0,\ldots,n+~1\}$ satisfying $u_{k+1}
\rightarrow_{D} u_{k}$ for $k\in \{0,\ldots i_0-1\}$, and $u_k
\rightarrow_{D} u_{k+1}$ for $k\in \{i_0,\ldots,n\}$.
Furthermore, if $j\geq i_0$, $u_j$ is proximate to some $u_{\sigma(j)}$ with
$\sigma(j)\leq i_0-1$.
\end{itemize}
\end{lem}

\begin{proof}
The first assertion of (a) is just Lemma 3.2 of \cite{FS2}. Now,
if $\u,\v\in ch_{\G}(\q,\p)$, either $\u\in ch_{\G}(\q,\v)$ or
$\v\in ch_{\G}(\q,\u)$; in any case, either $u$ is infinitely near
to $v$ or viceversa, and hence $u$ and $v$ cannot belong to
different branches of $D$. Now, we prove (b). First of all, note
that for any $i\in \{0,\ldots,n\}$ either $u_i$ is proximate to
$u_{i+1}$ or viceversa (cf. \ref{w=1extrem}). By (a), $u_1$ is
necessarily infinitely near to $u_0$ and so, proximate to it. If
each $u_{i+1}$ is proximate to $u_i$, the first claim is obvious
by taking $i_0=n+1$. Assume that there exists some $i\in
\{1,\ldots,n\}$ such that $u_i$ is proximate to $u_{i+1}$, and
take $i_0$ to be minimal with this property. We claim that
$u_{k+1} \rightarrow_{D} u_{k}$ for $k\in \{0,\ldots i_0-1\}$, and
$u_k \rightarrow_{D} u_{k+1}$ for $k\in \{i_0,\ldots,n\}$.
To show this, assume that there exists some $j\geq i_0+1$ such
that $u_{j+1}\rightarrow_{D} u_j$ and take $j_0$ to be minimal.
Then, both $u_{j_0-1}$ and $u_{j_0+1}$ are proximate to $u_{j_0}$
and, since they are adjacent to it, they are maximal among the
vertices of $D$ proximate to $u_{j_0}$. However, by (a) they are
in the same branch of $D$, so they must be equal, which is
impossible. Note that $u_{i_0}$ is the maximal point in $D$ among
the vertices belonging to $ch_{\G}(\q,\p)$. By (a) we know that
every $u_j$, $j\geq i_0$ is infinitely near to $q$. Write
$u_{\sigma(j)}$ for the maximal vertex in $D$ among the vertices
belonging to $ch_{\G}(\q,\u_{i_0-1})$ such that $u_j$ is
infinitely near to it. By (a) applied to
$ch_{\G}(\u_{\sigma(j)},\u_j)$ and the maximality of
$u_{\sigma(j)}$, necessarily $u_{\sigma(j)+1}$ is infinitely near
to $u_j$ and, because $u_{\sigma(j)+1}$ is proximate to
$u_{\sigma(j)}$, so is $u_j$. This completes the proof.
\end{proof}

The \emph{resolution graph} of a sandwiched singularity $\OO$ is
the weighted dual graph of the exceptional divisor of the minimal
resolution of $\OO$. These graphs are called \emph{sandwiched
graphs} and they are characterized as the weighted subgraphs of
some non-singular graph containing no vertices of weight $1$ (see
\cite{Spivak} Proposition II.1.11; cf. forthcoming \ref{rem_G}).
In particular, the graph obtained from a sandwiched graph by
removing an end is still a sandwiched graph.
\begin{rem} \label{rem_G}
If $\D$ is an Enriques diagram for $\OO$ and $\G^{0}_{D}$ is the
weighted subgraph of $\G_D$ comprising only the non-dicritical
vertices,
then $\G_{\OO}$ equals one of the connected components of
$\G^{0}_{D}$, whose vertices (and their corresponding vertices in
$D$) will be called \emph{non-dicritical vertices relative} to
$\OO$.
\end{rem}


\section{Sandwiched singularities and their Enriques diagrams} \label{section 2}

In Remark \ref{rem_G} we have observed that if $\D$ is an Enriques
diagram for a sandwiched surface singularity $\OO$, then its dual
graph $\G_D$ contains the resolution graph $\G_{\OO}$ as a
weighted subgraph.
Given any Enriques diagram $D$, the following proposition shows that
this combinatorial condition is sufficient to infer a result of geometrical
nature: suitable marking maps $\rho$ can be chosen so that $(D,
\rho)$ becomes an Enriques diagram for $\OO$.

\begin{pro} \label{graph-ideal}
Let $\OO$ be a sandwiched surface singularity and let $\D= (D,
\rho)$ be a marked Enriques diagram.
Then, $\D$ is an Enriques diagram for $\OO$ if and only if

\begin{enumerate}
\item the dual graph $\Gamma_{D}$ contains $\Gamma_{\OO}$ as a weighted subgraph.

\item $\rho (p) = 0$ if $\p \in \Gamma_{\OO}$; and $\rho (p) = +$ if $\p \in
\Gamma_{D}\setminus \Gamma_{\OO}$ and it is adjacent to some vertex of $\Gamma_{\OO}$.
\end{enumerate}
\end{pro}

\begin{proof}
The ``only if'' part follows from \ref{rem_G}, and, since
$\D=(D, \rho)$ is an Enriques diagram for $\OO$, the marking map
$\rho$ must satisfy the statement in order to assure that the
non-dicritical vertices of $D$ relative to $\OO$ correspond
exactly to the vertices of $\Gamma_{\OO}$.
For the ``if'' part, consider the reduced exceptional divisor
$E_{\OO}$ of the minimal resolution of $\OO$. Using plumbing (see
\cite{Lau72}, \cite{Spivak} Remark I.1.10), we can glue smooth
rational curves in order to obtain a configuration
$E_{\Gamma_{D}}$ containing $E_{\OO}$, having dual graph
$\Gamma_{D}$ and being contained on a smooth surface $S'$. Since
$\Gamma_{D}$ is a non-singular graph, by Castelnuovo's criterion
$E_{\Gamma_{D}}$ contracts to a non-singular point $O$
(\cite{Spivak} II.1.10) on a surface $S$. This contraction factors
into a composition of point blowing-ups (\cite{Lau71} theorem
5.7), so it determines a cluster with origin at $O$ and having
Enriques diagram $D$ (\cite{Cas00} 4.4). It is clear that $S'$ is
the surface obtained from $S$ by blowing up all points in $K$.
Consider a system of virtual multiplicities $\nu=\nu_K$ for $K$ in
such a way that all the points have positive excess except for those
corresponding to the vertices of the subgraph $\Gamma_{\OO}
\subset \Gamma_{D}$, which have excess 0 (such a system exists by
\cite{Cas00} 8.4.1).
Notice that no maximal point of $K$ corresponds to a vertex of $\Gamma_{\OO}$, since any
vertex of the resolution graph $\Gamma_{\OO}$ has weight strictly greater than $1$.
Hence ${\K}=(K,\nu)$ is a consistent cluster with positive virtual multiplicities. Let
$I=H_{\K}$ and $X=Bl_I(R)$. The morphism $f:S'\longrightarrow X$ given by the universal
property of blowing up is the minimal resolution of the singularities of $X$
(\cite{Spivak} II.1.4) and its exceptional components correspond to those points of
${\K}$ having excess 0. Therefore, $E_{\OO}$ is the exceptional divisor of $f$ and so,
the singularity on $X$ given by its contraction is isomorphic to ${\OO}$ (\cite{Lau71}
theorem 3.13).
\end{proof}

Given a sandwiched surface singularity $\OO$, there is not a
unique non-singular graph $\G$ containing $\G_{\OO}$ as a weighted
subgraph (see Example \ref{Ex1}). In fact, it is possible to
construct infinitely many non-singular graphs $\G$ containing a
given sandwiched graph $\G_{\OO}$, as Example \ref{Ex2} shows.

\begin{exam} \label{Ex2}
If $I\in \I$ is an ideal for $\OO$ and
$X=Bl_I(R)$, by choosing any non-singular point in the exceptional locus of $X$, and
blowing up this point, we obtain a new surface $X'$ containing $\OO$, as well. This $X'$
is the blowing-up of a complete $\m_O$-primary ideal $J_1=II_1\subset R$, where $I_1$ has
codimension one in $I$ (Theorem 3.5 of \cite{FS1}), and the dual graph $\Gamma_{J_1}$ contains a vertex more
than $\Gamma_I$.
In this way, an infinite chain of ideals in $\I$
\begin{equation*}
  \ldots \subset J_n \subset\ldots\subset J_1 \subset I \subset R
\end{equation*}
for $\OO$ can be constructed, and each $\G_{J_n}$ contains $\G_{\OO}$ as a weighted subgraph.
Moreover, for any $n$ the Enriques diagram $\D_{J_{n-1}}$ is a marked subdiagram of
$\D_{J_{n}}$.
\end{exam}

\begin{lem}\label{l.contraction}
Let $\D$ be an Enriques diagram for $\OO$. Consider the set $C$ of
non-dicritical vertices of $\D$ relative to $\OO$.
\begin{itemize}
\item[(a)] There is a tree structure on $C$ induced by the natural ordering $\leq$ of $D$.
\item[(b)] For any $p, \, q \in C$ define $p \rightarrow_{C} q$ if and only
if $p \rightarrow_{D} q$. Then,$\rightarrow_{C}$ is a proximity,
which turns $C$ into an Enriques diagram.
\end{itemize}
\end{lem}

\begin{proof}
To exhibit the tree structure of $C$, we will prove that
\begin{enumerate}
  \item there is a unique minimal element of $C$ by $\leq$, which is taken as the root of $C$;
  \item for any $p \in C$, its immediate predecessor in $C$ is the
  maximal element of $\{ q \in C: \ q< p \}$.
\end{enumerate}
Suppose that $p$ and $q$ are two different minimal vertices in
$C$, and write $w$ for the maximal vertex in $D(p)\cap D(q)$ (this
is, the maximal vertex which both $q$ and $p$ are infinitely near
or equal to). Then, as $\G_D$ contains no loops,
\[ch_{\G_D}(\q,\p)=ch_{\G_D}(\q,\w)\cup ch_{\G_D}(\w,\p),\]
and, by the connectivity of $\G_C$, we infer that $w\in C$,
contradicting the minimality of $q$ and $p$. We denote by $O_C$ the
minimal vertex of $C$, which is set as the root of $C$. On the
other hand, if $p\in C$, $p\neq O_C$, the vertices in $D(p)$ are
totally ordered by the natural ordering $\leq$ of $D$. Hence,
there exists a unique immediate predecessor of $p$, which is the
maximal element of $\{q\in C \mid q<p\}$, and this proves (a).

Now, to prove (b), we show that $\rightarrow_C$ defined as above
is a proximity relation for $C$. Since the root is the minimal
vertex of $C$, it is clear that it is proximate to no other vertex
of $C$. If $p\neq O_C$, its immediate predecessor $q_0$ in $C$ is
the maximal element of $\{q\in C \mid q<p\}$; hence $q_0<p$ and
$q_0\in C$. Then (b) of \ref{chain} says that $p$ is proximate to
some vertex $\w$ of $ch_{\G_D}(\q_0,\p)$, and (a) of \ref{chain}
says that $w$ is infinitely near or equal to $q_0$. Since
$ch_{\G_D}(\q_0,\p)\subset \G_C$, this leads to contradiction.
This proves the first condition of the proximity (see \S 1.3). The
conditions 2 and 3 are clearly satisfied.
\end{proof}

An Enriques diagram $C$ obtained as in Lemma \ref{l.contraction}
will be called \emph{a contraction for} $\OO$ (or \emph{for}
$\G_{\OO}$) \emph{associated with} $\D$. Reciprocally, we will
also say that $\D$ is \emph{associated with the contraction} $C$.

\begin{rem}\label{}
Notice that, by virtue of \ref{graph-ideal}, any Enriques diagram
(respectively, any contraction) for $\OO$ is in fact an Enriques
diagram (respectively, a contraction) for any sandwiched surface
singularity whose resolution graph is $\G_{\OO}$.
A contraction for $\OO$ may also be regarded as an enrichment of
the resolution graph $\G_{\OO}$ by some proximities between their
vertices, these proximities being compatible with the weights of
$\G_{\OO}$ in the sense of Lemma \ref{l.ineq-weights} below.
\end{rem}

\begin{exam} \label{Ex1}
Figure \ref{EDs} provides three distinct Enriques diagrams for the
same sandwiched singularity: they are not apparently related,
namely one is not a subdiagram of the other, as was the case in
Example \ref{Ex2}. Notice that $\D_1$ and $\D_3$ give rise to the
same contraction, which is the Enriques subdiagram of $D_1$
comprising the black dots.

\begin{figure}
\psfrag{A}{$\D_1$}\psfrag{B}{$\D_2$}\psfrag{C}{$\D_3$}
\begin{center}

\includegraphics[scale=1.25]{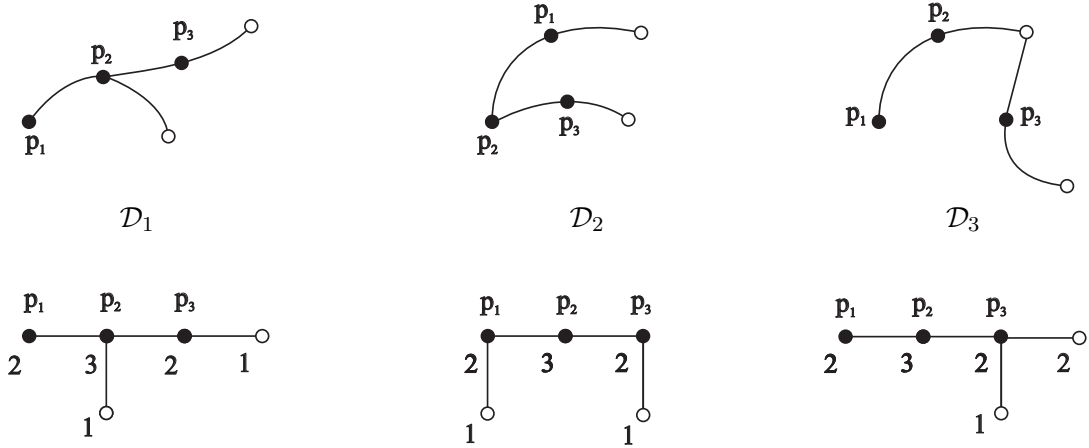}
\end{center}
\caption{\label{EDs}Three different marked Enriques diagrams for the same singularity
$\OO$ and their corresponding dual graphs. Dicritical vertices are represented with white
dots. }
\end{figure}
\end{exam}

\begin{rem}\label{r.contraction}
If $\D$ is an Enriques diagram for $\OO$ and $C$ is the associated
contraction, $C$ is not, in general, an Enriques subdiagram of $D$
(see Enriques diagram $\D_3$ of Figure \ref{EDs}). In particular,
if $I\in \I$ is an ideal for $\OO$ with Enriques diagram $\D$, the
set of points of $\K=BP(I)$ corresponding to the vertices of $C$
does not constitute, in general, a subcluster of $\K$.
\end{rem}

\begin{lem}\label{l.ineq-weights}
Let $C$ be a contraction for $\G_{\OO}$ associated with an
Enriques diagram $\D$.
Then for any vertex $p \in C$,
\begin{equation}\label{}
\omega_{\G_C}(\p) \leq \omega_{\Gamma_{D}}(\p) =
\omega_{\Gamma_{\OO}}(\p) \, ,
\end{equation}
and the inequality is strict at the extremal vertices of $C$. In
particular, $\G_C$ is not a weighted subgraph of $\G_{D}$.
\end{lem}

\begin{proof}
The inequality comes from the definition of contraction, since the
number of vertices proximate to $p$ in $C$ is less or equal than
in $D$. If $p$ is extremal in $C$, $\omega_{\G_C}(\p)=1$, while
$\omega_{\G_D}(\p)>1$, since $p$ is a non-dicritical vertex of
$\D$ and hence necessarily non-extremal in $D$ (see
\ref{w=1extrem}). The last assertion follows by considering the
weights at the extremal vertices of $C$.
\end{proof}

In \cite{Spivak} Corollary II.1.14, Spivakovsky introduced a type
of birational projection into a plane that could be achieved for
any sandwiched singularity. Namely he showed that, once a
sandwiched surface singularity $\OO$ is fixed, an ideal $I\in \I$
can be chosen in such a way that:
\begin{itemize}
 \item[\textbf{(i)}] $\OO$ is the only singularity of $X=Bl_I(R)$;
 \item[\textbf{(ii)}] the strict transform (by the minimal resolution of $X$)
 of any exceptional component of $\pi_{I}^{-1}(O)$ is a curve of the first
 kind, that is, the strict transform by $f$ (see diagram \ref{CommDiagr})
 of any  $L_p$ with $p \in BP(I)_{+}$
 has self-intersection equal to $-1$.
\end{itemize}
An ideal satisfying the above conditions (i) and (ii) (cf.
\cite{Mohring} 2.3) will be called an \emph{S-ideal for} $\OO$. A
marked Enriques diagram associated with an S-ideal for $\OO$ will
be called an \emph{S-Enriques diagram for} $\OO$.
The following result describes what the equisingularity classes of
S-Enriques diagrams look like:

\begin{lem}\label{idspivak}
An ideal $I\in \textbf{I}_R$ is an S-ideal if and only if the dicritical vertices of
$\D_I$ are free and extremal.
\end{lem}

\begin{proof}
Write $\D$ for the Enriques diagram of $\K=BP(I)$. First of all,
note that $X=Bl_I(R)$ has only one singularity $\OO$ if and only
if any non-dicritical vertex of $\G_D$ belongs to $\G_{\OO}$. Let
$p\in \K_+$ and assume that there exists some $q\in K$ infinitely
near to $p$. We may assume that $q$ is an immediate successor of
$p$. Then, $\omega_{\G_D}(\p)=r_D(p)+1\geq 2$ against condition
$(ii)$. Therefore, $p$ must be maximal in $K$. Now, assume that
$p$ is satellite, proximate to $u_1$ and $u_2$. Then, $\p\in
ch_{\G_D}(\u_1,\u_2)$. Necessarily, $u_1$ and $u_2$ are not
dicritical points of $K$ and thus, $\u_1,\u_2\in \G_{\OO}$. If
follows that $\p\in \G_{\OO}$ against the assumption $p\in \K_+$.

Conversely, if the dicritical vertices of $\D$ are free and
extremal, the union of the non-dicritical vertices of $\G_D$ is
connected and hence $X$ has only one singularity. Moreover, as
above, the self-intersection of the strict transform on $S_K$ of
any component $L_p$ with $p\in \K_+$ is $-1$.
\end{proof}

A contraction $C$ associated with an S-Enriques diagram will be
called an \emph{S-contraction}.
Contrary to what happened for general contractions (recall
\ref{r.contraction}), an S-contraction $C$ is a subdiagram
of its associated S-Enriques diagram $\D$;
furthermore, any S-contraction is associated with a unique
S-Enriques diagram:

\begin{pro}\label{P.S-contraction}
If $\D$ is an S-Enriques diagram for $\OO$, then an S-contraction
$C$ associated with $\D$ satisfies:
\begin{enumerate}
  \item[(a)] $C$ is a subdiagram of $D$;
  \item[(b)] $D$ can be recovered from $C$ by adding at any vertex $p \in C$ a number
  of $\omega_{\Gamma_{\OO}}(\p) - \omega_{\G_C}(\p)$ free
successors; $\D =(D, \rho)$ is recovered by defining the marking
map $\rho$ as $\rho (p)= 0$ if $p \in C$, and $\rho(p)=+$
otherwise.
\end{enumerate}
\end{pro}

\begin{proof}
By virtue of \ref{idspivak}, any dicritical point of $D$ is free
and extremal. Therefore, for any $p\in C$, there are exactly
$\omega_{\Gamma_{\OO}}(\p) - \omega_{\G_C}(\p)$ of these vertices
in the first neighborhood of $p$. This gives both claims.
\end{proof}

Next result is a sort of converse of Lemma \ref{l.ineq-weights}:

\begin{cor}\label{lem_2.2}
Let $C$ be an Enriques diagram and assume that the dual graph $\G_C$
equals $\G_{\OO}$ and satisfies $\omega_{\G_C}(\p) \leq
\omega_{\Gamma_{\OO}}(\p)$ at each vertex $p \in C$.
Then $C$ is an S-contraction for $\G_{\OO}$.
\end{cor}

\begin{proof}
To any $p\in C$ add exactly $\omega_{\Gamma_{\OO}}(\p) -
\omega_{\G_C}(\p)$ free successors to obtain from $C$ a marked
Enriques diagram $\D$ whose dicritical vertices are the extremal
ones. By construction the dicritical vertices are also free,
$\G_{\OO}$ is a weighted subgraph of $\G_{D}$. Then invoking
\ref{graph-ideal} and \ref{idspivak} we are done.
\end{proof}

Let us end this section by showing that the family of contractions
for a sandwiched surface singularity equals the family of
S-contractions:

\begin{pro} \label{GeneralitySpivIdeals}
Any contraction for $\G_{\OO}$ is an S-contraction.
\end{pro}

\begin{proof}
Let $\D$ be an Enriques diagram for $\OO$ and let $C$ be a
contraction for $\G_{\OO}$ associated with $\D$. Define a new
Enriques diagram $D'$ from $C$ by adding to each vertex $p\in C$
as many free successors as the number of vertices in $D \setminus
C$ that are proximate to $p$. Taking $\rho(p)=0$ if $p\in C$ and
$\rho(p)=+$ otherwise, the resulting Enriques diagram $\D'=(D',
\rho)$ for $\OO$ is an S-Enriques diagram associated with $C$ (see
\ref{graph-ideal}).
\end{proof}


\section{Contractions for a sandwiched surface singularity}

In this section we describe all the contractions for a given
sandwiched surface singularity $\OO$. Observe that, by virtue of
\ref{P.S-contraction} and \ref{GeneralitySpivIdeals}, this is
equivalent to listing all the equisingularity
classes of the S-ideals for $\OO$.

Suppose that the resolution graph $\Gamma _{\OO}$ has $n$ vertices
and that $\v$ is an end of $\G_{\OO}= \G_{n}$. The weighted graph
obtained by removing $\v$ is again a sandwiched graph and will be
denoted by $\Gamma_{n-1}$ (see \S \ref{ED-dualG}).
We want to detail a procedure to obtain all the contractions for
$\G_{n}$ from the contractions for $\Gamma_{n-1}$. Then, by
induction on $n$, the whole list of contractions for a given
sandwiched singularity will be inferred just from its resolution
graph.

The first result of this section describes how the vertex $v$
looks like in any contraction:
\begin{lem}\label{lem_2.1}
The vertex $v$ (corresponding to the end $\v$ of $\G_{\OO}$) in
any contraction $C$ for $\OO$ is either the root or free.
Furthermore, if $v$ is not the root of $C$, then either $v$ is
extremal or $v$ has a unique successor, which is satellite of $v$.
\end{lem}

\begin{proof}
Assume that $v$ is satellite, proximate to the vertices $u_1$ and
$u_2$ in $C$, and suppose that $u_2$ is proximate to $u_1$. We
will show that $\v\in ch_{\G_C}(\u_1,\u_2)$, thus contradicting
that $\v$ is an end, and proving the first claim.
Consider the Enriques subdiagram $C(v)$ of $C$ comprising all the
points preceding or equal to $v$. In particular, $u_1,u_2$ are
both in $C(v)$, and $v$ is maximal among the points of $C(v)$
proximate to $u_1$, and also to $u_2$. Hence, as vertices of
$\Gamma_{C(v)}$, $\v$ is adjacent to both $\u_1$ and $\u_2$ and
so, $\v\in ch_{\Gamma_{C(v)}}(\u_1,\u_2)$.
Now, the rest of vertices of $C \setminus C(v)$ all lie after some
vertex of $C(v)$, giving rise to blowing-ups of extra points.
The combinatorial effect of these blowing-ups is translated in the dual graph by the
elementary modifications introduced in I.1.5 of \cite{Spivak}, those of the first kind
corresponding to the blowing-ups of free points while those of the second kind to the
blowing-ups of satellite points. From their definition, it is immediate that these
modifications respect the property of being in the chain determined by two vertices
already in the graph.

For the second claim, let $\q$ be the only vertex in $\Gamma
_{\OO}$ to which $\v$ is adjacent, and assume that $v$ is not the
root of $C$.
We distinguish two cases. The first one is when $v$ is maximal
among the points in $C$ proximate to $q$. Since $v$ is an end,
there are no vertices in $C$ proximate to $v$ and $v$ is an
extremal vertex of $C$.
The second case is when $q$ is maximal among the points in $C$
proximate to $v$. Since $v$ is and end, $v$ must have a unique
successor, say $w$, preceding $q$. Denote by $u$ the immediate
predecessor of $v$. Since $v$ is free and $\v$ is adjacent only to
$\q$, $v$ cannot be the last point in $C$ proximate to $u$. Hence
$w$ must be proximate to $u$ and so $w$ is satellite of $v$.
This completes the proof.
\end{proof}

The following result shows how to construct a contraction for $\G_n$ from a
contraction for $\G_{n-1}$. Moreover, any contraction for $\G_n$
can be obtained in this way.

\begin{teo} \label{rules}
Suppose $\G_{n-1}$ is obtained from a sandwiched graph $\G_n$ by
removing an end $\v$.
Let $\u$ be the unique vertex in $\G_n$ to which $\v$ is adjacent,
and let $C'$ be a contraction for $\G_{n-1}$.
Define a new Enriques diagram by taking $C = C' \cup \{ v \}$ and
adding to $C'$ extra proximities relating $v$ according to one of
the following rules:
\begin{enumerate}
  \item If
$\omega_{\G_{C'}}(\u) < \omega_{\Gamma_{n}}(\u)$, set $v$ in $C$
as a free successor of $u$.
  \item If $u=q_r\rightarrow_{C'} q_{r-1} \rightarrow_{C'} \ldots \rightarrow_{C'} q_1$ are free
vertices in $C'$ with $1 \leq r< \omega_{\Gamma_{n}}(\v)$ and
$q_1\rightarrow_{C'} q_0$, then set $q_i \rightarrow_{C} v $ for
all $i \in \{1, \ldots , r \}$, and
either set $q_0 \rightarrow_{C} v$, provided $q_0$ is the root of
$C'$ and $r< \omega_{\Gamma_{n}}(\v)-1$ (with $v$ becoming the root of
$C$),
or set $v \rightarrow_{C} q_0$, provided
$\omega_{\G_{n-1}}(\q_{0}) < \omega_{\Gamma_{n}}(\q_{0})$.
\end{enumerate}
Then $C$ is a contraction for $\Gamma_{n}$. Moreover, any
contraction $C$ for $\G_n$ can be constructed from some
contraction $C'$ for $\G_{n-1}$ as above.
\end{teo}

\begin{proof}
Clearly $C$ defined as above satisfies $\G_C=\G_{n}$ and
$\omega_{\G_C}(\p) \leq \omega_{\Gamma_{n}}(\p)$ at each vertex $p
\in C$. Thus, invoking \ref{lem_2.2} and
\ref{GeneralitySpivIdeals}, the first claim follows.
By virtue of \ref{lem_2.1}, the vertex $v$ of $C$ corresponding to $\v$ is either free or
the root of $C$.
Since $\v$ is adjacent to $\u$ in $\Gamma _{n}$, there are only
two possibilities for their corresponding vertices in $C$:
\begin{description}
  \item[Case 1:] $v$ is maximal among the vertices in
  $C$ proximate to $u$.
  Hence $v$ cannot be the root of $C$, and by \ref{lem_2.1} $v$
  is free.
  Moreover, $v$ is an extremal vertex of $C$:
  otherwise, \ref{lem_2.1} implies that $v$ has a unique successor,
  which is satellite of $v$ and thus proximate to $u$,
  contradicting the maximality of $v$ among the vertices proximate to $u$.
  Therefore the set of vertices of $C'=
  C \smallsetminus \{ v \}$ is connected and has a tree structure.
  By considering the restriction of the proximity of $C$ to $C'$,
  $C'$ becomes an Enriques diagram.
  Clearly the graphs $\Gamma _{C'}$ and $\Gamma _{n-1}$ are equal
  (disregarding weights).
  Observe that $r_{C'}(u)=r_C(u)-1$ and that $r_{C'}(q)=r_C(q)-1$ if $q\in C$, $q\neq u$. Applying \ref{w=1extrem},
 we have
\begin{equation*}
  \omega _{\Gamma _{C'}} (\q) = \omega _{\Gamma _{C}} (\q) \leq \omega _{\Gamma _{n}} (\q)
  = \omega _{\Gamma _{n-1}} (\q) \, ,
\end{equation*}
  for any $q \in C' \smallsetminus \{ u \}$ and
\begin{equation*}
  \omega _{\Gamma _{C'}} (\u) = \omega _{\Gamma_{C} } (\u)-1 < \omega _{\Gamma _{C}} (\u) \leq \omega _{\Gamma _{n}} (\u)
  = \omega _{\Gamma _{n-1}} (\u) \, .
\end{equation*}
  Thus, invoking \ref{lem_2.2} and \ref{GeneralitySpivIdeals}, $C'$ is a contraction for $\G_{n-1}$.
  Finally, notice that the Enriques diagram
  $C$ is obtained from $C'$ by the
  procedure of the first rule of the statement,
  and we are done in this case.
  \item[Case 2:] $u$ is maximal among the vertices in
  $C$ proximate to $v$. Let $p_1, \ldots, p_j=u$ be
  the vertices in $C$ preceding or equal to $u$ which
  are proximate to $v$.
First of all, we  define on the set of vertices of $C'= C
\smallsetminus \{ v \}$ a tree structure.
We distinguish two cases:

\textbf{2.1.} If $v$ is the root of $C$, then take $p_1$ as the
root of $C'$, and for any $q \in C' \smallsetminus \{ p_1 \}$ declare
that $p$ is the immediate predecessor of $q$ in $C'$ if and only
if $p$ is the immediate predecessor of $q$ in $C$.

\textbf{2.2.} Otherwise, take the root of $C$ as the root of $C'$;
for any $q \in C' \smallsetminus \{ p_1 \}$ declare that
  $p$ is the immediate predecessor of $q$ in $C'$ if and only if $p$ is the immediate
  predecessor of $q$ in $C$; declare that the immediate predecessor of $p_1$ in $C'$
  is the immediate predecessor $p_0$ of $v$ in $C$.

Restrict the proximity of $C$ to $C'$, namely, for any $q,q'\in
C'$ set $q \rightarrow_{C'} q'$ if and only if $q \rightarrow_C
q'$. Let us check that it satisfies the properties 1 to 3 of a
proximity (see \S \ref{ED-dualG}). In the first case (where $v$ is
the root of $C$) these properties are clearly satisfied. In the
second case, the only condition that must be checked is property 1
for the vertex $p_1$, namely, that $p_1 \rightarrow_{C} p_0$.
Since $p_1$ is a successor of $v$ in $C$, by
  \ref{lem_2.1} we infer that $p_1$ is satellite of $v$. Thus $p_1$ is satellite in $C$: $p_1$ is
  proximate to its immediate predecessor in $C$, which is $v$, and to some point of $C$, say $p$;
  moreover, $v$ must be proximate to $p$, as well. Since $v$ is proximate to $p_0$, we infer that $p=p_0$ and
$p_1 \rightarrow_{C} p_0$, as desired.
  Therefore $C'$ is an Enriques diagram, whose dual
  graph $\Gamma _{C'}$ equals $\Gamma _{n-1}$ disregarding weights.
  Observe that $r_{C'}(p_0)=
  r_{C}(p_0)-1$ and $r_{C'}(q)= r_{C}(q)$ if $q\in D$, $q\neq p_0$. Applying \ref{w=1extrem},
 we have
\begin{equation*}
  \omega _{\Gamma _{C'}} (\q) = \omega _{\Gamma _{C}} (\q) \leq \omega _{\Gamma _{n}} (\q)
  = \omega _{\Gamma _{n-1}} (\q) \, ,
\end{equation*}
  for any $q \in C' \smallsetminus \{ p_0 \}$ and
\begin{equation*}
  \omega _{\Gamma _{C'}} (\p_0) = \omega _{\Gamma _{C}} (\p_0)-1 < \omega _{\Gamma _{C}} (\p_0)
   \leq \omega _{\Gamma _{n}} (\p_0)
  = \omega _{\Gamma _{n-1}} (\p_0) \, .
\end{equation*}
  Thus, invoking \ref{lem_2.2} and \ref{GeneralitySpivIdeals}, $C'$ is a contraction for
  $\G_{n-1}$.

  Finally, it remains to show that the Enriques diagram
  $C$ may be obtained from $C'$ by the
  procedure of the second rule of the statement.
  Indeed, according to the proximity defined in $C'$,
  notice first that $\{ u=p_j \rightarrow_{C'} \cdots
  \rightarrow_{C'} p_1 \}$ is a chain of free
  vertices in $C'$ preceding or equal to $u$,
  and that $p_1$ is the root of $C'$ if and only
  if $v$ is the root of $C$. On the other hand, recall that the proximity
  relations in $C$ involving the vertex $v$ are
\begin{equation*}
  p_{1} \rightarrow_{C} v , \ \ldots ,
  p_j \rightarrow_{C} v ,
\end{equation*}
and the further proximity relation $v \rightarrow_{C} p_0$ must be added in case $p_1$ is
not the root of $C'$.
This is exactly what performs the operation of the second rule,
and we are done.
\end{description}
\end{proof}

The whole list of contractions for $\G_{\OO}$ are obtained by
applying recursively the rules of \ref{rules}. Let us just sketch
the main steps of an implementation of this procedure. Each
step of this procedure adds a new vertex keeping the proximities
already defined. The idea is that the bigger the weights of
$\G_{\OO}$ are, the more Enriques diagrams for $\OO$ can be found.
\vspace{2mm}

\textbf{Step 1.} Choose any vertex of $\G_{\OO}$, say $\p_1$, and
take $\G_1=(\p_1,\omega_{\G_{\OO}}(\p_1))$ and
$C_1={\bullet}_{p_1}$. \vspace{2mm}

\textbf{Step i.} Assume that $\G_{i-1}$, $C_{i-1}$ have been
obtained, where $\G_{i-1}$ is a subgraph of $\G_{\OO}$. Choose a
vertex $\p_i$ adjacent to some $\q \in \G_{i-1}$. The graph $\G_i$
is obtained by adding $\p_i$ with weight $\omega_{\G_{\OO}}(\p_i)$
to $\G_{i-1}$, adjacent to $\q$; the new Enriques diagram $C_i$ is
obtained from $C_{i-1}$ by adding $p_i$ according to one of the
rules of \ref{rules}:

\begin{enumerate}
\item If
$\omega_{\G_{i-1}}(\q) < \omega_{\Gamma_{\OO}}(\q)$, $p_i$ can be
added to $C_{i-1}$ as a free successor of $q$,
$p_i\rightarrow_{C_i}q$.

\item If $q=q_r\rightarrow_{C_{i-1}} q_{r-1} \rightarrow_{C_{i-1}} \ldots \rightarrow_{C_{i-1}} q_1$ are free
vertices in $C_{i-1}$ with $1 \leq r< \omega_{\Gamma_{\OO}}(\p_i)$
and $q_1\rightarrow_{C_{i-1}} q_0$, then set $q_j
\rightarrow_{C_i} p_i $ for all $j \in \{1, \ldots , r \}$, and
either set $q_0 \rightarrow_{C_i} p_i$, provided $q_0$ is the root
of $C_{i-1}$ and $r< \omega_{\Gamma_{\OO}}(\p_i)-1$ (with $p_i$ becoming
the root of $C_i$),
or set $p_i \rightarrow_{C_i} q_0$, provided
$\omega_{\G_{i-1}}(\q_{0}) < \omega_{\Gamma_{\OO}}(\q_{0})$.
\end{enumerate}
The procedure stops at step $n$, the number of vertices of
$\G_{\OO}$. At this point, the obtained weighted graph $\G_n$
equals $\G_{\OO}$, and the Enriques diagram $C_n$ is just a
contraction for $\G_{\OO}$.

\begin{rem}
At any step of the procedure, there may be several choices to add
a fixed new vertex (for example, we may apply either rule 1 or 2
to add the vertex $p_i$ to $C_{i-1}$). In order to obtain the
whole list of all the contractions for $\G_{\OO}$, all these
possibilities must be performed. It might also happen that an
Enriques diagram $C_{i-1}$ to which the new vertex cannot be added
is reached. This means that no Enriques diagram for $\OO$ with the
subset of proximities of $C_{i-1}$ exists.
\end{rem}

\begin{rem}
Minimal singularities are rational surface singularities whose
fundamental cycle is reduced. They are characterized as those
sandwiched singularities that can be obtained by blowing up a
complete ideal all whose base points are free (see \cite{Mohring}
2.5; cf. \cite{Spivak}). As a consequence of our results, a
sandwiched surface singularity $\OO$ is minimal if and only if
there exists a contraction for $\OO$ that is obtained by applying
the first rule at each step of the above procedure.
\end{rem}

\begin{exam}\label{exSp}
Let $\OO$ be a singularity whose resolution graph is shown at the
bottom of Figure \ref{exSp}. By applying the procedure just
described, we obtain the whole list of contractions for $\OO$. The
S-Enriques diagrams shown in Figure \ref{llistatED} are obtained
by adding free successors to them as explained in (b) of
\ref{P.S-contraction}.

\begin{figure}
\begin{center}
\includegraphics[scale=1.]{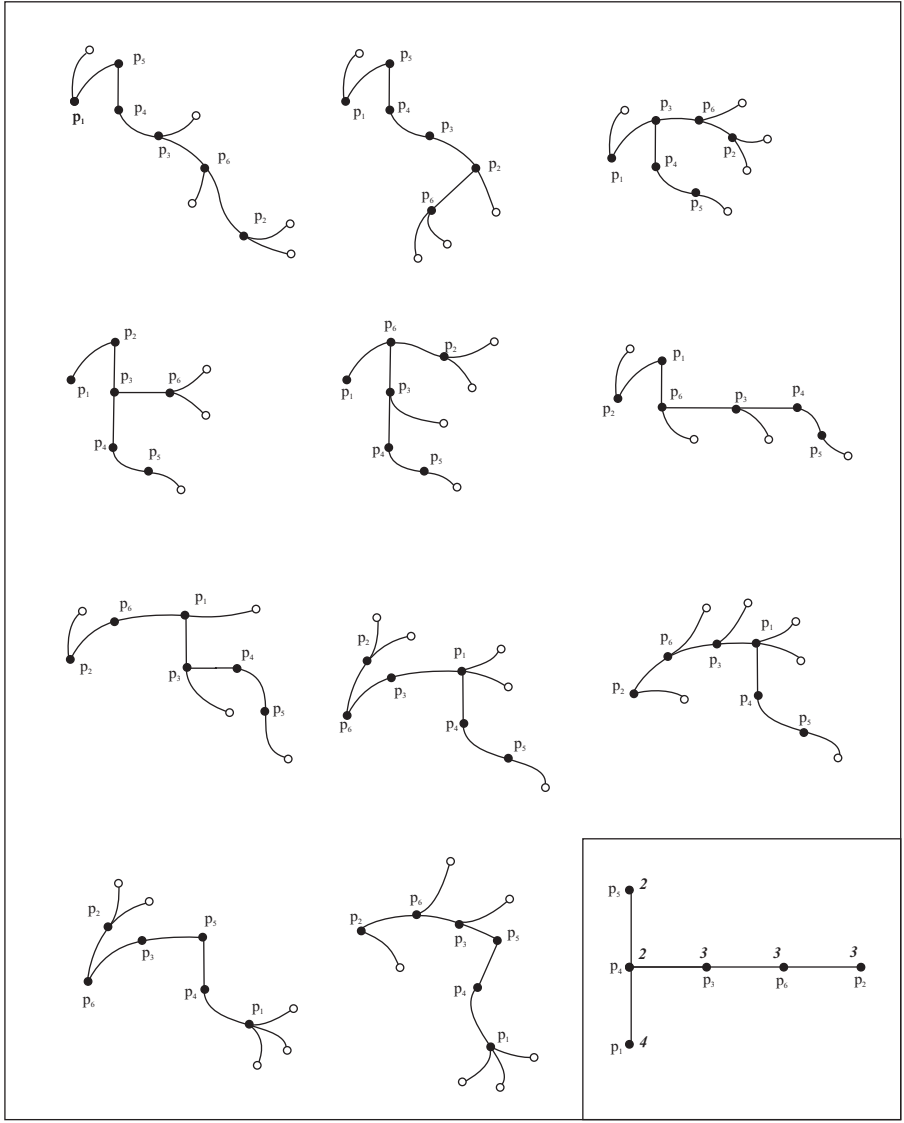}
\end{center}
\caption{ \label{llistatED} The complete list of S-Enriques diagrams for a singularity $\OO$
with resolution graph $\G_{\OO}$ in Example \ref{exSp}. The white-filled dots represent
the dicritical points, added to the contractions.}
\end{figure}

\end{exam}

\section{Equisingularity classes of the ideals for a sandwiched singularity}

In this section we address the problem of describing the
equisingularity classes of the ideals for a given sandwiched
surface singularity $\OO$, that is, of describing all the possible Enriques
diagrams for $\OO$.
The (finite) family of contractions for $\OO$ was inferred from
the resolution graph $\G_{\OO}$ of $\OO$ by the procedure
explained in the preceding section.
It remains to find out all the different Enriques diagrams
for $\OO$ giving rise to the same contraction (an infinite
family).
Here we will show how to complete contractions in order to
describe all the different Enriques diagrams for $\OO$,
thus solving completely the problem we are concerned with.

Given a contraction $C$ for $\OO$, our aim is to describe all the
Enriques diagrams for $\OO$ associated with $C$.
Consider the marked Enriques diagram $\C =(C, \rho_{C})$ with
$\rho_{C}(p)=0$ for any $p \in C$, and the number $\lambda_{C}=
\sum_{p \in C}(\omega_{\Gamma_{\OO}}(\p) - \omega_{\G_C}(\p))$.
By \ref{l.ineq-weights}, $\lambda_C >0$.
Let us describe a procedure to add vertices to $\C$ in order to
reach an Enriques diagram for $\OO$.
Write $\C_0 =(C_0, \rho_{0})=\C$.
For $1 \leq i \leq \lambda_C$, choose a vertex $p_i$ in $C$
such that $\omega_{\Gamma_{\OO}}(\p_i) >
\omega_{\G_{C_{i-1}}}(\p_i)$ and then define inductively
$\C_i =(C_i, \rho_{i})$ by taking $C_i= C_{i-1} \cup \{q_i
\}$ and ${\rho_{i}}_{|C_{i-1}} = \rho_{i-1}$, where the new
vertex $q_i$ is set as a successor of $p_i$ either
\begin{itemize}
  \item[A.] as a free successor of $p_i$, $q_i \rightarrow_{C_i} p_i$,
  and then set $\rho_i(q_i)=+$;
\end{itemize}
or, if there is some free successor $p'_i$ of $p_i$ in
$C_{i-1}$,

\begin{itemize}
  \item[B.] as a successor preceding $p'_i$, namely
  $q_i \rightarrow_{C_i} p_i$ and $p'_i \rightarrow_{C_i} q_i$ are the
  only proximities relating $q_i$, and then set $\rho_i (q_i)= +$
  in case $\rho_i (p'_i)= 0$ (otherwise $\rho_i (q_i)$ can be chosen no matter $0$ or
  $+$).
\end{itemize}

Notice that at step $i$ the operation of type A may
always be performed, independently of the existence of a free
successor of $p_i$, which would offer the possibility to
choose also an operation of type B.
Observe that $\sum_{p \in C}(\omega_{\Gamma_{\OO}}(\p) -
\omega_{\G_{C_i}}(\p))= \lambda_{C}-i$. Thus the procedure
performs effectively the $\lambda_{C}$ steps.
Any of such marked Enriques diagram $\C_{\lambda_{C}}$, obtained
from $C$ by the above procedure, will be called an
\emph{extension} of the contraction $C$.
Clearly any extension of $C$ is an Enriques diagram for $\OO$
associated with $C$.

\begin{rem}
Notice that any extension of $C$ all whose vertices have been
added performing operation A at each step is an S-Enriques diagram
for $\OO$ (in fact, the unique S-Enriques diagram for $\OO$
associated with $C$).
\end{rem}

The set of all extensions of $C$ forms a family of Enriques
diagrams for $\OO$ associated with $C$ minimal in the following
sense:

\begin{teo}\label{final}
Any Enriques diagram $\D$ for a sandwiched singularity $\OO$
contains, as a marked subdiagram, an extension of the contraction
associated with $\D$.

Conversely, if a marked Enriques diagram $\D $ contains, as a
marked subdiagram, an extension $\E$ of some contraction $C$ for
$\OO$ and satisfies that any vertex of $D\setminus E$ is proximate
to no vertex of $C$, then $\D$ is an Enriques diagram for $\OO$.
\end{teo}

\begin{proof}
For the first assertion, we need to find a marked subdiagram
$\E$ of $\D$ which is an extension of $C$.
Take $F= \{ p \in D : p \mbox{ is proximate to some }q \in C\}$,
and define $\E=(E, \rho_{E})$, where $E= C\cup F$ and $\rho_{E}$
is the restriction of $\rho_{D}$ to $E$.
Notice that $E$ is a connected subtree of $D$ since, if $p$ is
proximate to some $q\in C$, then any vertex in $D(p)$ infinitely
near to $q$ is also proximate to $q$. Hence, $E$ together with the
proximities inherited from the proximity of $D$ is an Enriques
subdiagram of $D$. Furthermore, $\E$ is a marked subdiagram of
$\D$.

Moreover, by \ref{w=1extrem}, the cardinality of $F$ equals
$\lambda:=\lambda_C$. Denote the vertices of $F$ by $\{ p_1,
\ldots, p_{\lambda}\}$ so that $p_i$ is not infinitely near to
$p_j$ if $j>i$. Write $\E_{\lambda}:= \E$ and for $1\geq
i<\lambda$, define, recursively $\E_i$ as the marked Enriques
diagram obtained from $\E_{i+1}$ by deleting $p_i$ (and keeping
the restricted proximity and marking map; the successors of $p_i$
become successors of the immediate predecessor of $p_i$). Notice
that the $\E_i$ are the marked Enriques diagrams generated by the
procedure detailed above to reach $\E$, proving that $\E$ is an
extension of $C$, as wanted.

For the converse, let $\E=(E, \rho_{E})$ be an extension of a
contraction $C$ for $\OO$. Thus, by \ref{graph-ideal}, $\G_{E}
\supseteq \G_{\OO} $ as weighted graphs, $\rho_{E}(p)=0$ for any
$\p \in \G_{\OO}$ and $\rho_{E}(p)=+$ for any $\p \in \G_E
\setminus \G_{\OO}$ being adjacent to some vertex of $\G_{\OO}$.
If $\D=(D, \rho_{D})$ contains $\E$ as a marked subdiagram and
there any vertex of $D \setminus E$  is proximate to no vertex of
$C\subset E$, then $\G_{D} \supseteq \G_{\OO} $ as weighted
graphs, and $\rho_{D}$ satisfies the marking map hypothesis of
\ref{graph-ideal}, 2: $\rho_{D}(p)=0$ for any $\p \in \G_{\OO}$
and $\rho_{D}(p)=+$ for any $\p \in \G_D \setminus \G_{\OO}$
adjacent to some vertex of $\G_{\OO}$. Hence, applying
\ref{graph-ideal} to $\D$ we are done.
\end{proof}

We have already pointed out that sandwiched singularities are
normal birational extensions of the regular ring $R$. If $R\subset
\OO$ is such an extension, there exists a complete ideal $I\subset
R$ such that $\OO=R[I/a]_{N_Q}$, where $N_Q$ is a height two
maximal ideal in $R[I/a]$ containing $\mathfrak{m}_R$ (the maximal
ideal of $R$), and $a$ is a generic element of $I$ (see
\cite{HunSally}). $R$ is said to be \emph{maximally regular} in
$\OO$ if there is no other regular ring $R'$ such that
\[R\subsetneq R' \subset \OO.\] Write $\D_I$ for the marked Enriques
diagram of the base points of $I$. Let $\E$ and $C$ be the
extension and the contraction for $\OO$ associated with $\D_I$.
Then, by virtue of \ref{final}, $\D_I$ can be thought as being
constructed from $\E$ by adding new vertices which are infinitely
near to some dicritical vertex of $\E$ and not proximate to any
vertex of $C$, or preceding the root of $C$ (notice that in any
case, the proximities of $\E$, and hence also the proximities of
$C$, are preserved). Moreover, $R$ is maximally regular in $\OO$ if and
only if the root of $\D_I$ equals the root of $\E$, i.e. no
vertices have been added to $\E$ preceding the root.

\begin{exam}
Let $\OO$ be a sandwiched singularity whose resolution graph is
shown in the top left corner of Figure \ref{llistatED2}. The
contractions for $\OO$ are shown at the top of the figure, and below
each one of them, a complete list of the associated extensions is
drawn. Any Enriques diagram for $\OO$ contains one of these
extensions as a marked subdiagram.

\begin{figure}

\begin{center}
\includegraphics[scale=0.78]{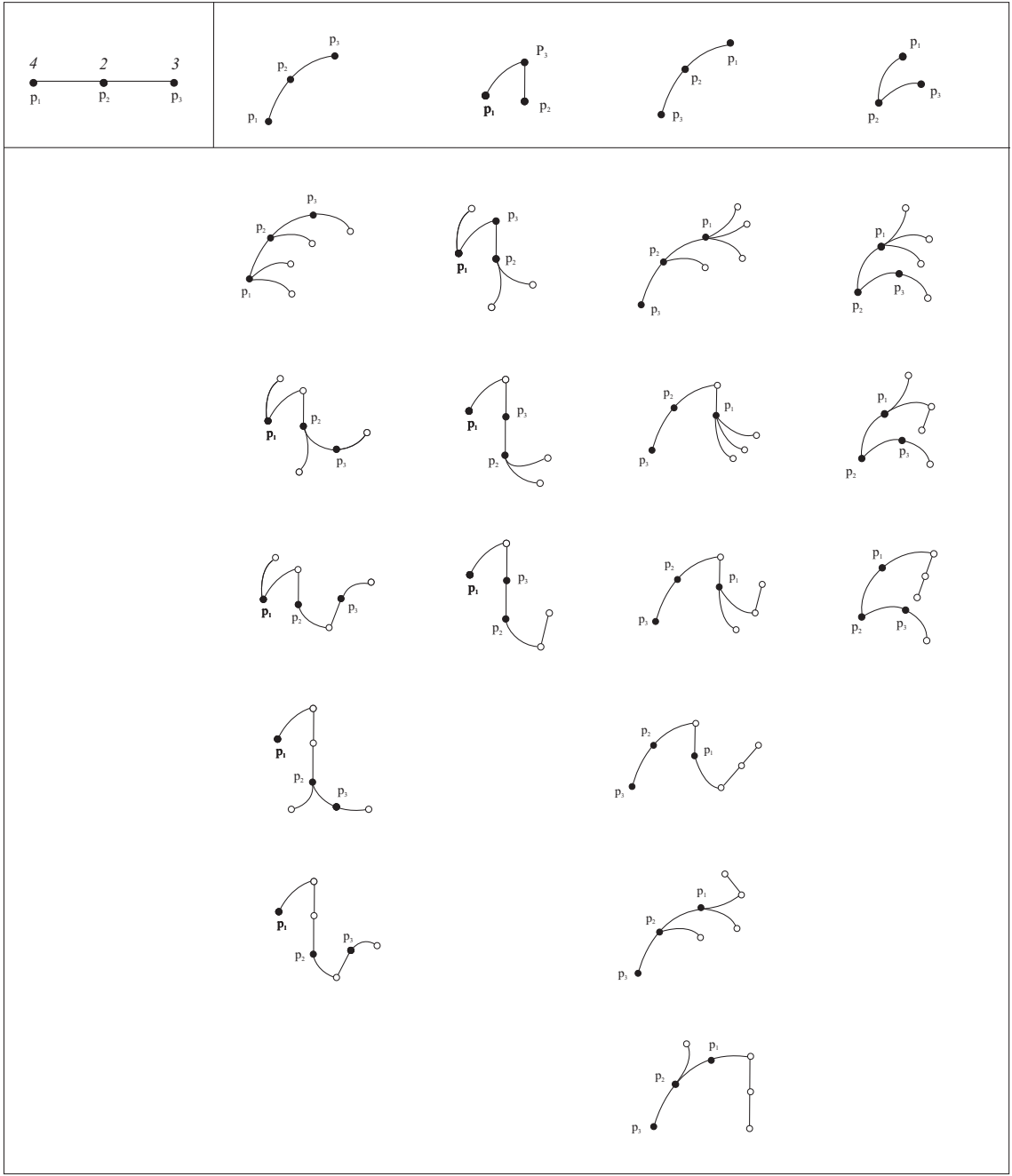}
\end{center}

\caption{\label{llistatED2} A complete list of extensions
for a sandwiched singularity whose resolution graph is drawn in the top left corner. At the top of the figure the
contractions of $\OO$ are shown. White-filled dots represent dicritical vertices.}
\end{figure}
\end{exam}

\nop  

\providecommand{\bysame}{\leavevmode\hbox to3em{\hrulefill}\thinspace}
\providecommand{\MR}{\relax\ifhmode\unskip\space\fi MR }
\providecommand{\MRhref}[2]{%
  \href{http://www.ams.org/mathscinet-getitem?mr=#1}{#2}
}
\providecommand{\href}[2]{#2}


\begin{thebibliography}{10}

\bibitem{Cas00}
E.~Casas-Alvero, \emph{Singularities of plane curves}, London Math. Soc.
  Lecture Notes Series, no. 276, Cambridge University Press, 2000.

\bibitem{EC85}
F.~Enriques and O.~Chisini, \emph{Lezioni sulla teoria geometrica delle
  equazioni e delle funzioni algebriche.}, Collana di Matematica [Mathematics
  Collection], vol.~5, Nicola Zanichelli Editore S.p.A., Bologna, 1985, Reprint
  of the 1915, 1918, 1924 and 1934 editions, in 2 volumes.

\bibitem{FS4}
J.~Fern{\'a}ndez-S{\'a}nchez,  \emph{On curves on sandwiched surface
  singularities}, (2006) arXiv: math/0701641.

\bibitem{FS1}
\bysame, \emph{On sandwiched singularities and complete ideals}, J. Pure Appl.
  Algebra \textbf{185} (2003), no.~1-3, 165--175.

\bibitem{FS2}
\bysame, \emph{Nash families of smooth arcs on a sandwiched singularity}, Math.
  Proc. Cambridge. Philos. Soc. \textbf{138} (2005), 117--128.

\bibitem{GSG92}
A.~Granja and T.~S\'anchez-Giralda, \emph{Enriques graphs of plane curves},
  Comm. Algebra \textbf{20} (1992), no.~2, 527--562.

\bibitem{HunSally}
C.~Huneke and J.~D. Sally, \emph{Birational extensions in dimension two and
  integrally closed ideals}, J. Algebra \textbf{115} (1988), no.~2, 481--500.

\bibitem{KP99}
S.~Kleiman and R.~Piene, \emph{Enumerating singular curves on surfaces}, Proc.
  Conference on {A}lgebraic {G}eometry: {H}irzebruch 70 (Warsaw 1998), vol.
  241, A.M.S. Contemp. Math., 1999, pp.~209--238.

\bibitem{Lau71}
H.~B. Laufer, \emph{Normal two-dimensional singularities}, Princeton University
  Press, Princeton, N.J., 1971, Annals of Mathematics Studies, No. 71.
  \MR{MR0320365 (47 \#8904)}

\bibitem{Lau72}
\bysame, \emph{On rational singularities}, Amer. J. Math. \textbf{94} (1972),
  597--608.

\bibitem{Lip69}
J.~Lipman, \emph{Rational singularities, with applications to algebraic
  surfaces and unique factorization}, Inst. Hautes \'Etudes Sci. Publ. Math.
  (1969), no.~36, 195--279.

\bibitem{Mohring}
K.~M{\"o}hring, \emph{On sandwiched singularities}, Ph.D. thesis, November
  2003.

\bibitem{Reguera97}
A.J. Reguera-L\'{o}pez, \emph{Curves and proximity on rational surface
  singularities}, J. Pure Appl. Algebra \textbf{122} (1997), no.~1-2, 107--126.

\bibitem{Spivak}
M.~Spivakovsky, \emph{Sandwiched singularities and desingularization of
  surfaces by normalized {N}ash transformations}, Ann. of Math. (2)
  \textbf{131} (1990), no.~3, 411--491.

\end{thebibliography}
\end{document}